\documentclass[1p,times]{elsarticle} 

\makeatletter
\def\ps@pprintTitle{%
  \let\@oddhead\@empty
  \let\@evenhead\@empty
  \def\@oddfoot{}%
\let\@evenfoot\@oddfoot}
\makeatother

\usepackage{enumerate}
\usepackage{amssymb,amsmath,amsfonts}
\usepackage{natbib}

\usepackage[most]{tcolorbox}

\definecolor{darkred}{rgb}{0.7,0,0}
\definecolor{darkgreen}{rgb}{0,0.7,0}
\newcommand{\Red}[1]{{\color{darkred}#1}}
\newcommand{\Green}[1]{{\color{darkgreen}#1}}

\usepackage{xcolor}
\newcommand{\commentout}[1]{}

\usepackage{booktabs}

\usepackage{graphicx}
\usepackage{amssymb}
\usepackage{hanging}
\usepackage{url}

\title{A Bayesian approach to calibrating hydrogen flame kinetics using many
experiments and parameters}
\author[1]{John Bell}
\author[1]{Marcus Day}
\author[2]{Jonathan Goodman}
\author[3]{Ray Grout\corref{cor}}
\author[4]{Matthias Morzfeld}
\address[1]{Computational Research Division, Lawrence Berkeley National Lab, Berkeley, CA 94720, USA}
\address[2]{Courant Institute of Mathematical Sciences, New York University, New York, NY 10012, USA}
\address[3]{Computational Science Center, National Renewable Energy Laboratory, Golden, CO 80401, USA}
\address[4]{Department of Mathematics, University of Arizona, Tucson, AZ, 85721, USA}
\cortext[cor]{Corresponding author}
\ead{ray.grout@nrel.gov}

\begin{document}
\begin{frontmatter}

\begin{abstract}
First-principles Markov Chain Monte Carlo
sampling is used to investigate uncertainty quantification and uncertainty propagation in
parameters describing hydrogen kinetics.
Specifically, we sample the posterior distribution of thirty-one parameters
focusing on the H$_2$O$_2$ and HO$_2$ reactions
resulting from conditioning on ninety-one experiments.
Established literature values are used for the remaining parameters in the mechanism.
The samples are computed using an affine invariant sampler starting with broad, noninformative priors.
Autocorrelation analysis shows that O(1M) samples are sufficient to obtain a reasonable sampling of the posterior.
The resulting distribution identifies strong positive and negative correlations
and several non-Gaussian characteristics.
Using samples drawn from the posterior, we 
investigate the impact of parameter uncertainty on the prediction of two more complex flames: a 
2D premixed flame kernel and the ignition of a hydrogen jet issuing into a heated chamber.
The former represents a combustion regime similar to the target
experiments used to calibrate the mechanism and the latter 
represents a different combustion regime.  For the premixed flame, the net amount of
product after a given time interval has a standard deviation of less than 2\% whereas the
standard deviation of the ignition time for the jet is more than 10\%.  
The samples used for these studies are posted online.
These results indicate the degree to which parameters consistent with the target experiments
constrain predicted behavior in different combustion regimes.
This process provides a framework for both identifying reactions for further study from candidate mechanisms as well as combining uncertainty quantification and propagation to, ultimately, tie uncertainty in laboratory flame experiments to uncertainty in end-use numerical predictions of more complicated scenarios.

\end{abstract}


\end{frontmatter}

\section{Introduction}
Assembling combustion kinetic mechanisms is an arduous task, in which
computation and simulation have increasing importance.
A mechanism is specified by giving the functional form of the rate
equations and a ``calibration'' that determines numerical values for 
the unknown parameters.
Traditionally, parameters in the rate expressions (the ``calibration'' in
our terminology) are determined from direct measurements designed to isolate
specific reactions, and from theory. These experiements are typically difficult to design since most observable quantities are the result of many reaction steps and first-princples calculation of elementrary rates from quantum mechanics principles is an active research area. 
The resulting mechanism is then tested
against simplified laboratory flame experiments, with sensitivity analysis
used to identify how predictions depend on specific parameters.
The sensitivity analysis, combined with estimated error bars on parameters, can
be used to identify reactions that require further study.

A complementary approach is to augment direct experiments and theory with 0D (homogeneous ignition) 
and 1D (steady flame) experiments and integrate that data into the mechanism design and 
calibration process. 
In these experiments, the highly nonlinear parameter dependence can be determined only 
through simulation.  
The calibration is optimized to find the best fit to a set of flame experiments.
Prior knowledge about the reactions is incorporated into the operation as inequality
bounds on the parameters.
This approach was used to develop GRIMech \cite{grimech} and is currently used by a number
of groups in the chemical kinetics community \cite{FrenklachWR92, DavisJWE05, LiYWL15}.
To reduce the computational effort required for the optimization, most groups
have adopted a response surface approach.  In this approach, sensitivity analysis is used
to determine what parameters are most important for a given experiment and a
response surface model of the experiment is constructed based on that reduced set of parameters.

Here, we adopt the philosophical approach of Bayesian uncertainty quantification, see e.g.,
\cite{Stuart:2010} and 
many references there.
Rather than finding the best
fit to kinetic parameters, our goal is to characterize the distribution of parameters
that are consistent with a given set of experimental data.
This approach allows us to assess not only which parameters are important for a given
set of
experiments but also which parameters are not.  
Capturing the broader distribution of parameters consistent with the data provides a natural
framework for assessing the predictive capabilities of the resulting model.
Results from this type of study can be used to identify which reactions have the greatest impact on uncertainty and evaluate
the impact of additional experiments on reducing that uncertainty.

Our methodology uses first principles Bayesian uncertainty quantification followed by 
uncertainty propagation.
The uncertainty quantification methodology is similar to that used by Braman et al. \cite{BramanOR13} 
to quantify the distribution of parameters in several syngas mechanisms consistent
with a database of premixed laminar flames.
But there are several technical differences.
One difference is that we use broad uninformative priors.
Informative priors lead to posteriors that depend strongly on subjective judgements based 
on presupposed knowledge of the parameters.
When the posterior depends strongly on prior bounds, we believe this may indicate that the
new experiments partly contradict previous understanding and therefore may be given precedence.
Unfortunately, broad priors lead to much larger regions of parameter space to explore, which
puts stringent demands on the sampling methodology. We have used an affine invariant ensemble
sampler \cite{emcee_hammer} that can handle wide ranges in scales, strong correlations between
parameters and non-Gaussian behavior.
The present study would have been impossible without this type of methodology.

Here we will focus on hydrogen kinetics.
Although hydrogen kinetics have been studied extensively,
some aspects of the hydrogen oxygen mechanism are still uncertain.
The present study focuses on 31 calibration parameters related to pre-exponential factors
and third body coefficients 
in the $\mbox{H}\mbox{O}_2$ 
and $\mbox{H}_2\mbox{O}_2$ pathways that are thought to be significant at high pressure.
Given the complex interactions between reactions and parameters, we believe it is more informative 
to study the simultaneous impact of many parameters than to isolate them into smaller subgroups.
We note that our decision to use broad priors for this study disregards current knowledge of these
parameters.  Here our goal was to eliminate strong dependence of the posterior distribution on
the prior for the selected parameters and assess what can be learned about the parameters
solely from the given target experiments.
For all of the other parameters describing the mechanism, we used the accepted values
in the literature and exclude them from the sampling, effectively adopting a zero-variance prior
for these parameters.
As well as shedding light on a particularly thorny set of reactions, this is an examplar for how to treat a plausible situation in mechanim development. Often a number of reactions are well understood (e.g., through elementary rate experiments or quantum  calculations) while others are less well characterized. For those that are less well characterized it is useful to have a rigorous procedure to guide further study. The procedure adopted here can be used to identify parameters that the experiments are insensitive to as well as to understand if there are alternative combinations of parameter sets that are equally consistent with the data. Equipped with this knoweldge, the unimporant parameters can be left with very approximate values and new experiments can be designed that can discriminate between the pathways that are indistinguishable with existing data. 

The present study uses published data from 91 experiments, primarily laminar flame
speed measurements augmented with a small number of flow reactor and ignition delay experiments.
The Bayesian posterior is defined in terms of the mismatch between computational predictions
and measured exerimental data normalized by experimental error bars.
Each likelihood requires computation of steady flames and time-dependent 0D flow reactor and ignition delay experiments.
We use the {\tt PREMIX} code \cite{PREMIX98} for steady flames, and {\tt VODE} for the 0D 
and ignition delay simulations.
After a long MCMC run and an auto-correlation analysis of the output, we used the samples to represent the posterior distribution as a set of points. Each sample represents a collection of 31 parameters that is consistent with the 91 experiments.
 This representation of the posterior makes it easy to  examine projections into lower dimensions to facilitate visualization and also draw samples from the posterior by selecting points at random while avoiding any constraints on the form of the distribution. 

We also explored the potential impact of posterior parameter uncertainty using computational
uncertainty propagation. 
We chose two numerical experiments whose results would be sensitive to kinetics effects and
ran them for several hundred parameter sets that characterize the posterior.
We emphasize that the parameter sets for the different runs are different, but all of them
are consistent with the experiments.
The first experiment is propagation of a flame in a pre-mixed medium from a small 
crenelated kernel.
The resulting flame is thermodiffusively unstable, making it senstitive to kinetics; however,
it is in a similar combustion regime to many of the target experiments used for calibration.
For this case, although the flames show some differences in structure, the net fuel consumption is fairly insensitive
to the different choices of parameters.
The other experiments is the ignition of a jet of cold $\mbox{H}_2$ fuel into a bath of a hot mixture of $\mbox{O}_2$
and $\mbox{N}_2$.
This non-premixed ignition experiment represents a combustion regime that is not
well represented by the calibration experiments.
Not surprisingly, in this case
the predicted ignition time of the jet shows a strong dependence on the parameters,
indicating that for this type of problem there is still significant uncertainty in the predictions.

The results of our Bayesian calibration are posted online in the form of 576 sample calibrations.
These are the samples used for the uncertainty propagation experiments just described.
In the Bayesian philosophy, the posterior distribution is the calibration.
This set of samples seems to be the best representation of the posterior that we can offer.
We forego the traditional practice of publishing the best fit calibration and error bars because we
do not feel that this is an accurate representation of the posterior.
The posterior means and covariances are readily estimated from the samples.

The organization of the rest of the paper is as follows.
Section 2 contains a detailed description of the problem setup.
This includes a discussion of the kinetic mechanism and a discussion of our choice of parameters.
We also discuss the experiments that were selected for the study along with the associated
experimental error distributions. Finally we give a precise definition of the prior distributions.
In Section 3, the affine invariant MCMC sampler that we selected is motivated and described.
Section 4 presents an analysis of the MCMC output.
This includes an auto-correlation analysis that studies the quality of the resulting chain.
It also includes some statistical analysis of the samples to explore the posterior distribution.
Finally, in Section 5, we describe the two uncertainty propagation experiments and present the results.

\section{Problem setup}

\subsection{Hydrogen kinetics}
\label{sec:H2}

Hydrogen-oxygen kinetics have been the subject of extensive study,
both because of interest in hydrogen and hydrogen-enriched fuels and
because the H$_2$/O$_2$ mechanism forms an important submechanism in
hydrocarbon kinetic mechanisms.

Two primary datasources have led to
the current state-of-the-art hydrogen oxidation mechanisms. First,
elementary rate measurements, taken from carefully constructed experiments
sensitive to one or a small number of individual reaction rates are
used to infer individual rate parameters such as done by
\cite{MuellerYD98, AshmanH98}. Secondly, macro experiments where the
observable depends on the entire reaction mechanism have been used to
validate and understand reaction pathways through mechanisms created by assembling the elementary reactions. 
Such mechanisms have been developed and updated as individual rates have been refined, e.g. the sequence of work 
by Yetter et al. \cite{YetterDR91},  Mueller et al.  \cite{MuellerKYD99}. 

In addition to exploring how the various elementary reactions interact, 
comparison of macro observables to experimental measurements has been used for refinement of the mechanism in a `comprehensive' sense, such as undertaken by 
 by Li et al.  \cite{LiZKD04} and O'Conaire et
al. \cite{ConaireCSPW04}. The Li et al.  mechanism has been
updated by Burke et al. \cite{BurkeCDJ10, BurkeCJDK12} while the
O'Conaire et al. mechanism has been updated by Keromnes et
al. \cite{Keromnes_Many_Curran13}. A similar kinetic mechanism has
been developed by Konnov \cite{Konnov08}, also by assembling
elementary reactions from the literature.

Closely related developments have been made for kinetic models
involving hydrogen and carbon monoxide (syngas). Li et
al. \cite{LiZKCDS07} augmented their earlier H$_2$/O$_2$ mechanisms with
C$_1$/O$_2$ kinetics and then adjusted key rate constants identified by
sensitivity analysis to improve predictions of the macro
experiments. Davis et al.  \cite{DavisJWE05} placed somewhat heavier
emphasis on macro optimization and developed a H$_2$/CO mechanism
based on assembling a trial mechanism drawing significantly the
H$_2$/C$_1$/O$_2$ chemistry from GRI-Mech \cite{grimech} and then
optimizing a significant number of kinetics parameters against a
library of laminar flame, flow reactor and shock-tube experiments. You
et al. \cite{YouPF11} produced a viable hydrogen mechanism as a
byproduct of establishing the `data collaboration' method as part of a
workflow to automate mechanism generation. The mechanism was
subsequently the focus of an uncertainty analysis study performed by
Li et al. \cite{LiYWL15}.

We note that the literature cited above represents two significantly different approaches to
mechanism development.  One approach
relies on adjustment of individual reactions done in the context of insight and analysis of
reaction pathway fluxes
whereas the other places greater reliance on sensitivity studies and global optimization to
provide a mechanism that targets a specific set of macro experiments.
The Bayesian approach used here represents, in some sense, an intermediate approach that attempts to
develop a more general purpose mechanism and uses macro experiments to constrain parameters rather than
fit parameters.  Although not pursued here, detailed reaction path analysis can be
incorporated into the Bayesian framework, which would represent a more blended approach.

\subsection{Choice of parameters}

For the analysis here, we have started with the basic reaction model developed by
Burke et al. \cite{BurkeCDJ10, BurkeCJDK12}.
We have augmented that basic reaction set by incorporating two additional reactions
identified by Burke as X1 and X6. The baseline mechanism is summarized in Table \ref{tab:mech}.
For this baseline mechanism, we have chosen to vary the pre-exponential factors and third-body
coefficients for reactions involving HO$_2$ and H$_2$O$_2$.  These active parameters are indicated in
red and green in Table \ref{tab:mech}.  (The two green coefficients, the Ar and He third body coefficients,
were set equal to each other and represent a single parameter.)
These parameters play an important role
in high-pressure lean flames and are not as well understood as the reactions involving H, O and OH radicals.

As noted above, for the present study we will select broad priors so that the prior does little to constrain the
portion of parameter space that is explored. (We again acknowledge that this 
choice ignores much of the prior art.)
The prior is a product of one variable priors for each parameter.
This makes the parameters independent in the prior.
Figure \ref{fig:tp} shows strong dependences between variables in the posterior.
We first select a Gaussian for each parameter with a mean given by the value of the parameter
as specified by Burke et al.\ and a standard deviation equal to the mean. 
For reactions X1 and X6, we use the characterization of the reactions in You et al.\ \cite{YouPF11} to set the mean values.
There values are summarized in Table \ref{tab:mech}.
We then restrict these priors with upper and lower bounds to prohibit nonphysical parameter choices, such as
negative pre-exponential factors or extremely large values that lead to excessive failures of the software
used to evaluate the different experiments.  These upper and lower bounds are also summarized in the table.
Thus, the prior for each parameter is a ``truncated'' Gaussian.

\renewcommand{\baselinestretch}{1.1}
\begin{table}
\scriptsize
  \centering
   \caption{Arrhenius rate parameters for syngas combustion model (kinetics, and accompanying 
           thermodynamics and transport -- in \cite{CHEMKINIII96} format).  Parameter database taken from
           \cite{BurkeCJDK12}), except for the final two, which were taken from \cite{YouPF11}.
           Parameters in \Red{red} and \Green{green} are active for the present study;
           the green parameters are varied synchronously.
           The forward rate constant, $K_f = AT^{\beta}\exp{(-E_a/RT)}$. $^*$The number
           in parentheses is the exponent of 10, i.e., 2.65(16) = 2.65 $\times$ 10$^{16}$.
           For the active parameters lower and upper bounds (hidden constraints) are provided in the last column.}
  \label{tab:mech}

  \begin{tabular}{r @{\hspace{1em}} r c l @{\hspace{1em}} r @{\hspace{1em}} r @{\hspace{1em}} r @{\hspace{2em}} r}
\multicolumn{1}{c}{RID}  &
\multicolumn{3}{c}{Reaction} &
\multicolumn{1}{c}{A$^{*}$} &
\multicolumn{1}{c}{$\beta$} &
\multicolumn{1}{c}{$E_a$} &
Bounds L:U \\ \hline
R1  & H + O$_2$       &=& O + OH             & 1.04(14)        &  0      &15286 & \\
R2  & O + H$_2$       &=& H + OH             & 3.818(12)       &  0      & 7948 & \\
    & O + H$_2$       &=& H + OH             & 8.792(14)       &  0      &19170 & \\
R3  & OH + H$_2$      &=& H + H$_2$O         & 2.16(8)         &  1.51   & 3430 & \\
R4  & 2 OH            &=& O + H$_2$O         & 3.34(4)         &  2.42   &-1930 & \\
R5  & H$_2$ + M       &=& 2 H + M            & 4.577(19)       & -1.40   &104380 & \\
    &  \multicolumn{6}{@{\hspace{1.2em}} l}{{\em Third-body:}  H$_2$(2.5), H$_2$O(12), Ar(0), He(0)}& \\
    & H$_2$ + Ar      &=& 2 H + Ar           & 5.84(18)        & -1.10   &104380& \\
    & H$_2$ + He      &=& 2 H + He           & 5.84(18)        & -1.10   &104380& \\
R6  & 2 O + M         &=& O$_2$ + M          & 6.165(15)       & -0.5    &    0 & \\
    &  \multicolumn{6}{@{\hspace{1.2em}} l}{{\em Third-body:}  H$_2$(2.5), H$_2$O(12), Ar(0), He(0)}& \\
    & 2 O + Ar        &=& O$_2$ + Ar         & 1.886(13)       &  0      &-1788 & \\
    & 2 O + He        &=& O$_2$ + He         & 1.886(13)       &  0      &-1788 & \\
R7  & O + H + M       &=& OH + M             & 4.714(18)       & -1      &    0 & \\
    &  \multicolumn{6}{@{\hspace{1.2em}} l}{{\em Third-body:}  H$_2$(2.5), H$_2$O(12), Ar(0.75), He(0.75)}& \\
R8  & H$_2$O + M      &=& H + OH + M         & 6.064(27)       & -3.322  &120790 & \\
    &  \multicolumn{6}{@{\hspace{1.2em}} l}{{\em Third-body:}  H$_2$(3), H$_2$O(0), He(1.1), N$_2$(2), O$_2$(1.5)}& \\
    & H$_2$O + H$_2$O &=& H + OH + H$_2$O    & 1.006(26)       &  -2.44  &120180& \\
R9  & H + O$_2$ (+M)  &=& HO$_2$ (+M)        &  &  & & \\
    &  \multicolumn{3}{@{\hspace{1.2em}} l}{{\em high pressure, K$_{f\infty}$}} & \Red{4.65084(12)} & 0.44 & 0  & 2(12):1(13) \\
    &  \multicolumn{3}{@{\hspace{1.2em}} l}{{\em low pressure, K$_{f0}$}}      & \Red{6.366(20)} & -1.72 & 524.8  & 0:8(20)\\
    &  \multicolumn{6}{@{\hspace{1.2em}} l}{{\em TROE:} $F_c = 0.5$} & \\
    &  \multicolumn{3}{@{\hspace{1.2em}} l}{{\em Third-body:}} & H$_2$(\Red{2.0}) & & & 0:6\\
    &  \multicolumn{3}{@{\hspace{1.2em}} l}{\hspace{5em}} & H$_2$O(\Red{14}) & & & 0:28\\
    &  \multicolumn{3}{@{\hspace{1.2em}} l}{\hspace{5em}} & O$_2$(\Red{0.78}) & & & 0:3\\
    &  \multicolumn{3}{@{\hspace{1.2em}} l}{\hspace{5em}} & Ar(\Red{0.67}) & & & 0:3\\
    &  \multicolumn{3}{@{\hspace{1.2em}} l}{\hspace{5em}} & He(\Red{0.8}) & & & 0:3\\
R10  & HO$_2$ + H      &=& H$_2$ + O$_2$      & \Red{2.750(6)}        &  2.09   & -1451 & 1(6):5(16)\\
R11  & HO$_2$ + H      &=& 2 OH               & \Red{7.079(13)}       &  0      &  295 & 2(13):1(14)\\
R12  & HO$_2$ + O      &=& OH + O$_2$         & \Red{2.850(10)}        &  1      &-723.93 & 1(9):1(11)\\
R13  & HO$_2$ + OH     &=& O$_2$ + H$_2$O     & \Red{2.890(13)}        &  0      & -497 & 1(13):6(13)\\
R14  & 2 HO$_2$        &=& O$_2$ + H$_2$O$_2$ & \Red{4.200(14)}       &  0      &11982 & 1(14):2(15)\\
    & 2 HO$_2$        &=& O$_2$ + H$_2$O$_2$ & \Red{1.300(11)}       &  0      &-1630 & 5(10):4(11)\\
R15  & H$_2$O$_2$ (+M) &=& 2 OH (+M)       &  & & \\
    &  \multicolumn{3}{@{\hspace{1.2em}} l}{{\em high pressure, K$_{f\infty}$}} &  \Red{2.00(12)} & 0.9   & 48749 & 5(11):1(12)\\
    &  \multicolumn{3}{@{\hspace{1.2em}} l}{{\em low pressure, K$_{f0}$}}      &  \Red{2.49(24)} & -2.3  & 48749 & 1(23):1(25)\\
    &  \multicolumn{6}{@{\hspace{1.2em}} l}{{\em TROE:}
              $F_c=0.43$}& \\
    &  \multicolumn{3}{@{\hspace{1.2em}} l}{{\em Third-body:}} & H$_2$(\Red{3.7}) & & & 0:15\\
    &  \multicolumn{3}{@{\hspace{1.2em}} l}{\hspace{5em}} & H$_2$O(\Red{7.5})& & & 0:20\\
    &  \multicolumn{3}{@{\hspace{1.2em}} l}{\hspace{5em}} & H$_2$O$_2$(\Red{7.7})& & & 0:20\\
    &  \multicolumn{3}{@{\hspace{1.2em}} l}{\hspace{5em}} & O$_2$(\Red{1.2})& & & 0:5\\
    &  \multicolumn{3}{@{\hspace{1.2em}} l}{\hspace{5em}} & N$_2$(\Red{1.5})& & & 0:5\\
    &  \multicolumn{3}{@{\hspace{1.2em}} l}{\hspace{5em}} & He(\Red{0.65})& & & 0:4\\
R16  & H$_2$O$_2$ + H  &=& OH + H$_2$O        & \Red{2.410(13)}       &  0      & 3970 & 5(11):1(14)\\
R17  & H$_2$O$_2$ + H  &=& HO$_2$ + H$_2$     & \Red{4.820(13)}       &  0      & 7950 & 1(12):9(13)\\
R18  & H$_2$O$_2$ + O  &=& OH + HO$_2$        & \Red{9.550(6)}        &  2      & 3970 & 1(5):3(7)\\
R18  & H$_2$O$_2$ + OH &=& HO$_2$ + H$_2$O    & \Red{1.740(12)}       &  0      &  318 & 5(10:5(12)\\
    & H$_2$O$_2$ + OH &=& HO$_2$ + H$_2$O    & \Red{7.590(13)}       &  0      & 7270 & 4(12):4(14)\\
X1  & HO$_2$ + H      &=& H$_2$O + O         & \Red{3.97(12)}        &  0      &  671 & 1(12):9(12)\\
X6  & O + OH + M      &=& HO$_2$ + M         & \Red{8.000(15)}       &  0      &   0 & 2(15):2(16)\\
    &  \multicolumn{3}{@{\hspace{1.2em}} l}{{\em Third-body:}} & H$_2$(\Red{2})& & & 0:6\\
    &  \multicolumn{3}{@{\hspace{1.2em}} l}{\hspace{5em}} & H$_2$O(\Red{12}) & & & 0:35\\
    &  \multicolumn{3}{@{\hspace{1.2em}} l}{\hspace{5em}} & Ar(\Green{0.7}) & & & 0:3\\
    &  \multicolumn{3}{@{\hspace{1.2em}} l}{\hspace{5em}} & He(\Green{0.7}) & & & 0:3\\
  \end{tabular}
 \end{table}

\subsection{Choice of experiments}

We have selected 91 experiments to provide data for the calibration.  Of these experiments, 77 are
laminar premixed flame experiments across a broad range of stoichiometries and pressure.
Of these 77 premixed flames, 71 were taken from Burke et al.\ \cite{BurkeCJDK12}. The remaining 6 premixed flames
were taken from (Refs 27-31 in \cite{DavisJWE05}).
Data from the 77 premixed flames were augmented with measurements from 14 experiments used
by Davis, et al.\ for optimization of a syngas mechanism, including 6 flow reactors (Ref 13 in \cite{DavisJWE05})
and 8 ignition experiments (Refs 36-39 in \cite{DavisJWE05}).
The premixed flame cases were simulated using PREMIX \cite{PREMIX98};
for each case we extracted the propagation speed of the steady, unstrained flame.
The flow reactor and ignition delay experiments were simulated as point (0D) reactors at
constant pressure or volume, respectively, using the backward-difference integration method in VODE \cite{vode};
the data extracted from each of these simulations is detailed in the relevant references,
as are the measured values and experimental errors bars.
For each experiment, we have assumed that measurement errors can be characterized by a Gaussian distribution with
mean given by the experimental measurement and variance set by experimental error bars reported in the experiment.

\section{Sampling Methodology}  \label{sec:MCMC}

Our goal is to characterize the distribution of reaction parameters that is consistent
with the available data.  
We will adopt a Bayesian sampling approach that allows us to sample this distribution.
To make this precise, we define a
combustion model ${\cal M}(\boldsymbol{\theta})$
that computes a predicted model vector $\mathbf{z_\theta}$
given a vector $\boldsymbol{\theta}$ of active parameters.
Here the number of active parameters $n_{\theta} = 31$  and the output vector
is of size $n_z = 91$.
We assume that we are given measurements $d_k$ for the suite of experiments
and denote the standard deviation of the measurement error by $\sigma_k$.
Given a prior distribution, $p(\boldsymbol{\theta})$ of $\boldsymbol{\theta}$,
we can estimate the likelihood that $\mathbf{z_\theta}$ will match experimentally measured data.
Using Bayes' rule (see e.g.\ \cite{Stuart:2010})
\[
p(\boldsymbol{\theta} | \mathbf{d}) = \frac{1}{Z_f} p_{\boldsymbol{\theta}} (\boldsymbol{\theta}) \,
p(\mathbf{d} | \boldsymbol{\theta})
\]
where
$p(\mathbf{d} | \boldsymbol{\theta})$ is the distribution of predicted data given $\boldsymbol{\theta}$ and
$Z_f$ is normalization factor.
Here we assume that measurement uncertainty for experiment $k$ is a Gaussian with mean ${d}_k$ 
and variance $\sigma_k^2$ so that
\begin{equation}
     p(d|\boldsymbol{\theta}) =
\frac{1}{Z_L} \exp\!\left( - \sum_{k=1}^{n_z} \frac{\left( {z}_{\theta,k} - {d}_k \right)^2}{2\sigma_k^2}\right) \; .
\label{eq:post}
\end{equation}

If one assumes $p_{\boldsymbol{\theta}} (\boldsymbol{\theta})$ is
Gaussian and ${\cal M}$ is linear, then $p(\boldsymbol{\theta} | \mathbf{d})$ is also Gaussian, and so its posterior would be completely determined
by its mean and variance. The mean can be found by minimizing the quantity
$$
  F(\boldsymbol{\theta}) = - \log p(\boldsymbol{\theta} | \mathbf{d}) =
-\log p_{\boldsymbol{\theta}} (\boldsymbol{\theta}) - \log p(\mathbf{d} | \boldsymbol{\theta}) - \log(Z_f) \; .
$$
Define the minimizer, $\boldsymbol{\mu} =\text{arg}\, \min_{\boldsymbol{\theta}} F(\boldsymbol{\theta})$.
When the model is nonlinear, this {\em posterior mode} (also called MAP point, for maximum a-posterior point)
need not be the posterior mean.
There is no guarantee that there is only a single global posterior mode.
The optimization package {\tt MINPACK} was unable to reliably identify the posterior mode.

We have used Markov Chain Monte Carlo (MCMC) to produce samples of the posterior distribution (\ref{eq:post}).
Specifically, we use the emcee hammer software package of Foreman-Mackey
\cite{emcee_hammer}, which is based on
the stretch-move ensemble sampler of Goodman and Weare \cite{GoodmanWeare:2010}.
The stretch move ensemble sampler addresses several challenges in MCMC sampling.
First, the sampler is affine invariant, which makes the sampler
good at long narrow valleys and avoids the need to precondition the problem
by a change of sampling variables; such optimizations typically require hand-tuning.
Second, the sampler uses multiple walkers that can be divided into two groups each of which
can be evaluated in parallel.

\section{Sampling study}  \label{sec:AC}

This section has two parts.
First we describe the results of the MCMC process.
Auto-correlation studies show that many MCMC steps are needed for good sampling.
We explore possible reasons for this.
Next we visualize the posterior sample.

We ran the emcee hammer code with an ensemble size $L=64$ for approximately
$T=1.5 \times 10^4$ steps.
This produced $9.6 \times 10^5$ samples, many of which are highly correlated.
The initial ensemble is constructed by sampling each parameter from a Gaussian
distribution with the same mean as the prior but with the variance reduced by
a factor of 100.
The sampler has a single dimensionless tuning parameter, $a>1$, that sets the length 
scale of proposed moves as a multiple of distance between samples.
We took $a=1.3$ throughout. 
We used the red-black parallel version of the sampler, which allows for 
$32 = 64/2$ independent likelihood evaluations per re-sampling sweep.

\subsection{Autocorelation}

We assess the quality of the samples by studying the auto-correlation functions
of the individual variables in the calibration.
We refer the reader to \cite{Sokal:1989} for background on the role of auto-correlation
studies in analysis of MCMC output.
We write $\boldsymbol{\theta}_{i,k,t}$ for the value of parameter $i$ in walker $k$ after $t$
sampler sweeps through the ensemble.
Then $\boldsymbol{\theta}_{k,t}= (\boldsymbol{\theta}_{1,k,t}, \ldots, \boldsymbol{\theta}_{d,k,t})$ is the sample calibration 
corresponding to walker $k$ at sweep $t$.
For each $i$ and $k$, the numbers $\boldsymbol{\theta}_{i,k,t}$ form a time series with a 
theoretical auto-covariance function
\begin{equation}   \label{eq:ACdef}
        C_{i,s} = \mbox{cov}(\boldsymbol{\theta}_{i,k,t}, \boldsymbol{\theta}_{i,k,t+s}) \; .
\end{equation}
The theoretical lag $s$ auto-covariance is the same for each walker, $k$.
We estimate (\ref{eq:ACdef}) as follows.
We omit the first $T_b = T/2$ samples, $T_b$ being a burn-in time.
This unusually large burn-in time reflects the impact of using a tight cluster of initial sample points to 
initialize the sample and the long decay times of the 
auto-covariance function we observed.
For each walker and parameter, we estimate the empirical auto-covariance as
\[
       \widehat{C}_{i,k,s} 
          = \frac{1}{T-T_b - s} \sum_{t=T_b}^{T-s}\left( \boldsymbol{\theta}_{i,k,t} - \overline{\boldsymbol{\theta}}_{i,k}\right)
                        \left( \boldsymbol{\theta}_{i,k,t+s} - \overline{\boldsymbol{\theta}}_{i,k}\right) \; ,
\]
with $\overline{\boldsymbol{\theta}}_{i,k}$ being the sample mean.
Our overall estimate of $C_{i,s}$ is found from these by averaging over walkers in the ensemble
\[
        \widehat{C}_{i,s} = \frac{1}{L}  \sum_{k=1}^L \widehat{C}_{i,k,s} \; .
\]
The auto-correlations are the auto-covariances normalized by the lag zero covariance:
\begin{equation}  \label{eq:rho}
        \widehat{\rho}_{i,s} = \frac{\widehat{C}_{i,s}}{\widehat{C}_{i,0}} \; .
\end{equation}
It is a peculiar feature of the emcee hammer algorithm, which is not shared by other samplers,
that the auto-correlation functions of different parameters are similar.

\begin{figure}
  \centering
  \includegraphics[width=.5\textwidth]{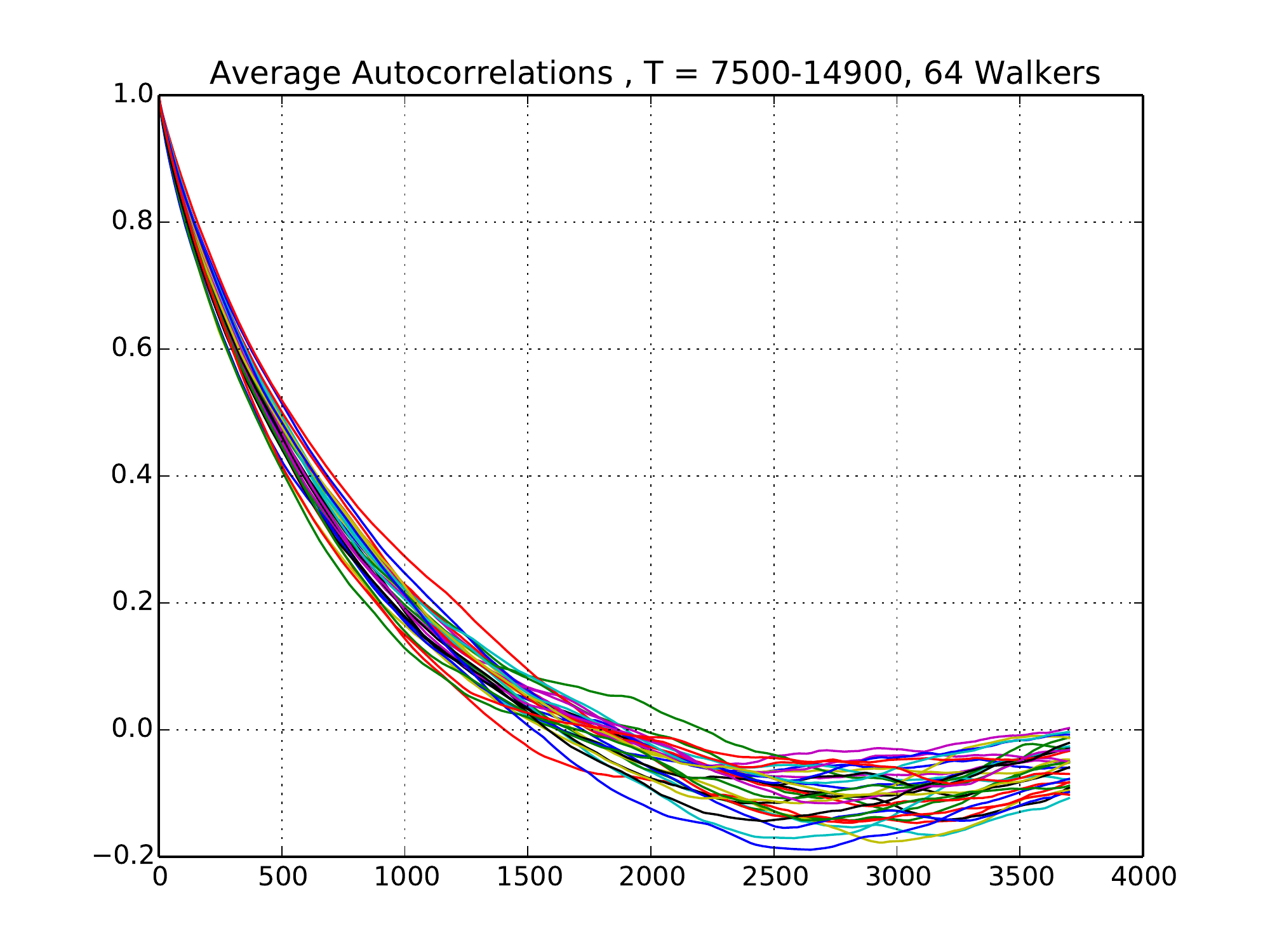}
  \caption{Average autocorrelation of 64 walkers for each parameter}
  \label{fig:AC}
\end{figure}

This is evident in Figure \ref{fig:AC}, which plots $\widehat{\rho}_{i,s}$ for all $d = 31$
parameters.
It shows that the auto-correlations decay on a time scale of a thousand sweeps.
This suggests that our run of $T=15 (10^3)$ sweeps produces a modest number of
effectively independent ensembles.
In principle, we should give a quantitative estimate of the auto-correlation time.
But we judged that our run was too short (measured in auto-correlation times) to make
such an estimate reliable.
The samples that we  used for uncertainty propagation were extracted from the second half of the run.

The long auto-correlation times indicate that this posterior distribution is particularly hard 
to sample. 
This makes a striking contrast with other experience with the sampler, which is extensive.  See \cite{emcee_hammer} for
references to applications.
The triangle plots below do not indicate that the distribution is strongly multi-modal, but
such effects may only be visible in higher dimensions.

\subsection{Posterior}  \label{sec:post}

As we already stated, we represent the posterior as a point cloud consisting of a number 
of sample calibrations. 
It is a challenging research area to visualize the properties of a point cloud in 
31 dimensions.
The triangle plot (described below) may be used to explore single parameter distributions 
and pairwise correlations.
We emphasize that many more subtle correlations and relationships between parameters are possible.
Uncertainty propagation is likely to be the only reliable way to determine how parameter 
uncertainties and correlations effect a target experiment.

A triangle plot is a collection of histograms and scatterplots. 
It is organized as a 2D array, with the plots indexed by parameter pairs, $(i,j)$,
with $i$ on the horizontal axis and $j$ on the vertical.
For each $j < i$ we show the histogram of pairs $(X_{i,k,t}, X_{j,k,t})$ over all $k$ and $t$.
\begin{figure}[tph]
  \centering
  \includegraphics[width=.9\textwidth]{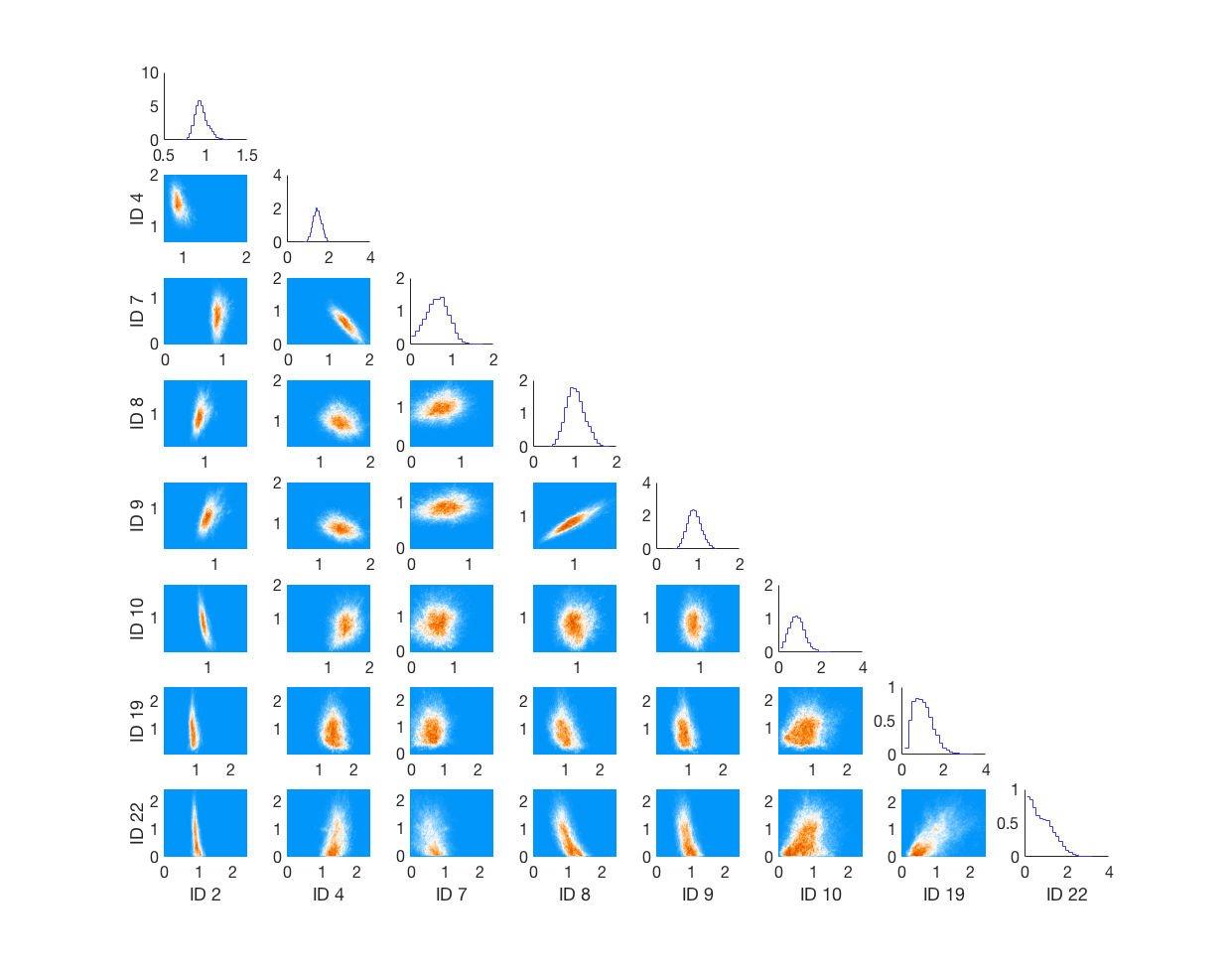}
  \caption{Triangle plot for subset of the parameters sampled ($i,j 
  \in (2,4,7,8,9,10,19,22)$ ). Axis are normalized by prior means.}
  \label{fig:tp}
\end{figure}
On the diagonal positions, $(i,i)$, we show a histogram of the parameter $X_{i,k,t}$.
For 31 parameters, the individual plots in the triangle plot shown on a printed page
in this format would be too small to be very useful.
Figure \ref{fig:tp} is a sparse triangle plot, drawn using eight of the variables chosen to be representative of the qualitative features of the full triangle plot. The parameters chosen are shown in Table~\ref{tab:tpparams}.

\begin{table}[h!]
  \label{tab:tpparams}
  \caption{Subset of active parameters appearing in triangle plot in Figure~\ref{fig:tp}}
  \centering
\begin{tabular}{l l}
ID & Role \\
\hline
2  & Pre-exponential factor in low pressure part of R9 \\
4  & Third body efficiency of H$_2$O in R9 \\
7  & Third body efficiency of He in R9 \\
8  & Pre-exponential factor in R10 \\
9  & Pre-exponential factor in R11 \\
10 & Pre-exponential factor in R12 \\
19 & Third body efficiency of N$_2$ in R15 \\
22 & Pre-exponential factor in R16 \\
\end{tabular}
\end{table}

The last two parameters have highly skewed distributions and
the bottom two rows of Figure \ref{fig:tp} in particular display strongly non-Gaussian behavior.
The bottom row is indicative of a parameter where the posterior apparently has non-negligible likelihood right up to a parameter bound.
The second row from the bottom on the other hand displays a non-Gaussian character that drops to zero before encountering the parameter bound.
With the exception of several parameters discussed below,
the posterior is contained within the parameter bounds. That is, the entire posterior is captured within the bounds of the priors, indicating that the priors is sufficiently broad
that they do not effect overall results.
Other interesting qualitative behavior includes shapes that indicate positively correlated (9,8); negatively correlated (7,4); inverse (10,2), (22,2), (22,8);  Gaussian (10,9), (10,9); and complex relationships (19,10), (22,10). Many of the pairs do not seem to have strong correlations. Not included in Figure~\ref{fig:tp} but of interest because of lack of dependence is the high pressure pre-exponential factor in  R9 and, to a lesser degree, the pre-exponential factor in reaction X1. In both cases the posterior is broad and tends towards uniformity. This is interesting as it suggests these parameters are not constrained by the data: the experiments chosen are simply not sensitive to these parameters. 

Some of observed behavior can be traced back to physically meaningful relationships. For example, the positive correlation in $(9,8)$ is between the pre-exponential factors in competing reactions R10 and R11: \\
\begin{tabular}{r @{\hspace{1em}} r c l }
 \bf{R10}  & HO$_2$ + H      &$\Leftrightarrow$& H$_2$ + O$_2$  \\
\bf{R11}  & HO$_2$ + H      &$\Leftrightarrow$& 2 OH            \\
\end{tabular} \\
Previous work (\cite{BurkeDJ11}) noted that the competition between these two reactions contributes substantially to the pressure dependent behavior of the mechanism. The positive correlation in the distribution of the posterior reflects a constraint on the relative flux through these competing pathways. 

The negative correlation between parameters $(7,4)$ involves the third body efficiency of He and H$_2$O in reaction R9: \\
\begin{tabular}{r @{\hspace{1em}} r c l }
\bf{R9}  & H + O$_2$ (+M)  &$\Leftrightarrow$& HO$_2$ (+M);  \\
\end{tabular} \\
where the net rate of this reaction can be maintained by equivalently by either third body.
A further observation is that for several of the parameters (e.g., ID 4, 10, 19, 22) the mode of the posterior has shifted substantially from the prior mean. The behavior in the bottom row, where the mode is at the lower bound (recall from Table~\ref{tab:mech} that the lower bound is two orders of magnitude less than the prior mean), suggests that the mechanism might be more consistent with the data if this reaction was removed.

\section{Uncertainty Propagation Results} \label{sec:UP}

In the previous section, we showed that there is a wide range of mechanism parameters that
are consistent with our target experiments.  An important question to ask is how much do 
the resulting uncertainties in the mechanism influence predictions for other experiments.
To quantify the resulting uncertainty, we select samples from the posterior distribution 
and use those samples to perform simulations of additonal experiments.  Statistical analysis
of the resulting simulations provides an assessment of the uncertainty that does not rely on
any assumptions about the distribution of the uncertainty.

\begin{figure}
  \centering
  \includegraphics[width=.3\textwidth]{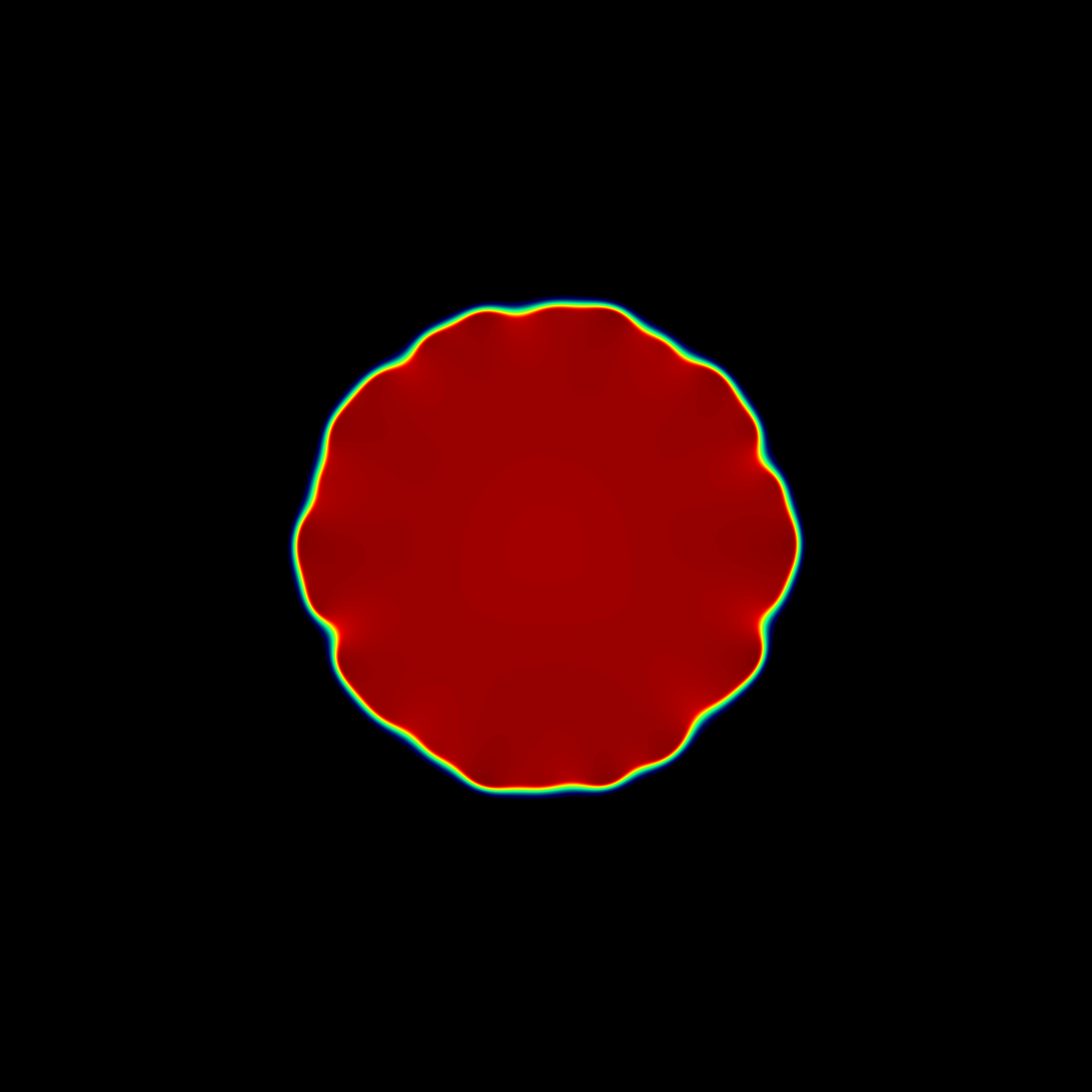}
  \includegraphics[width=.3\textwidth]{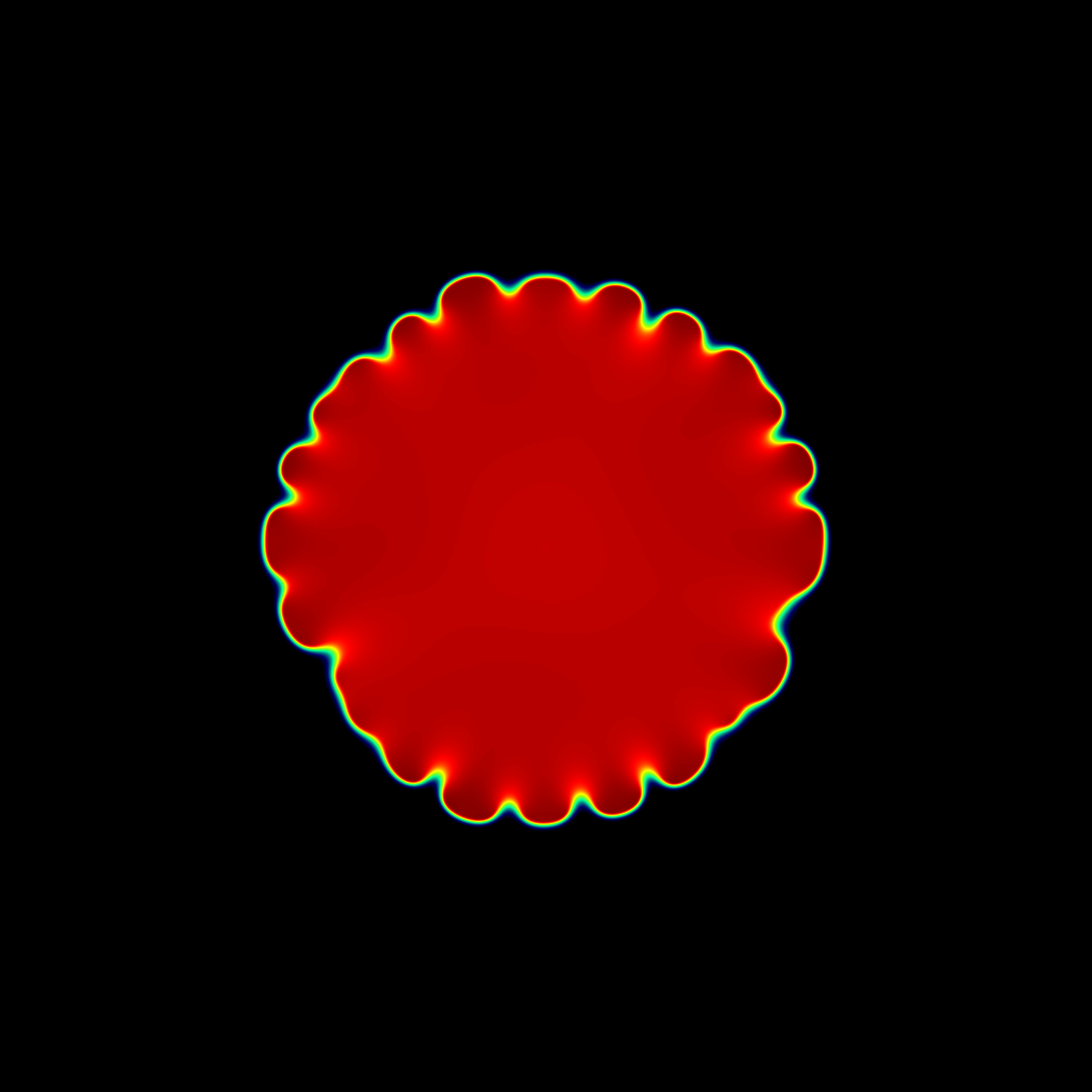}
  \includegraphics[width=.3\textwidth]{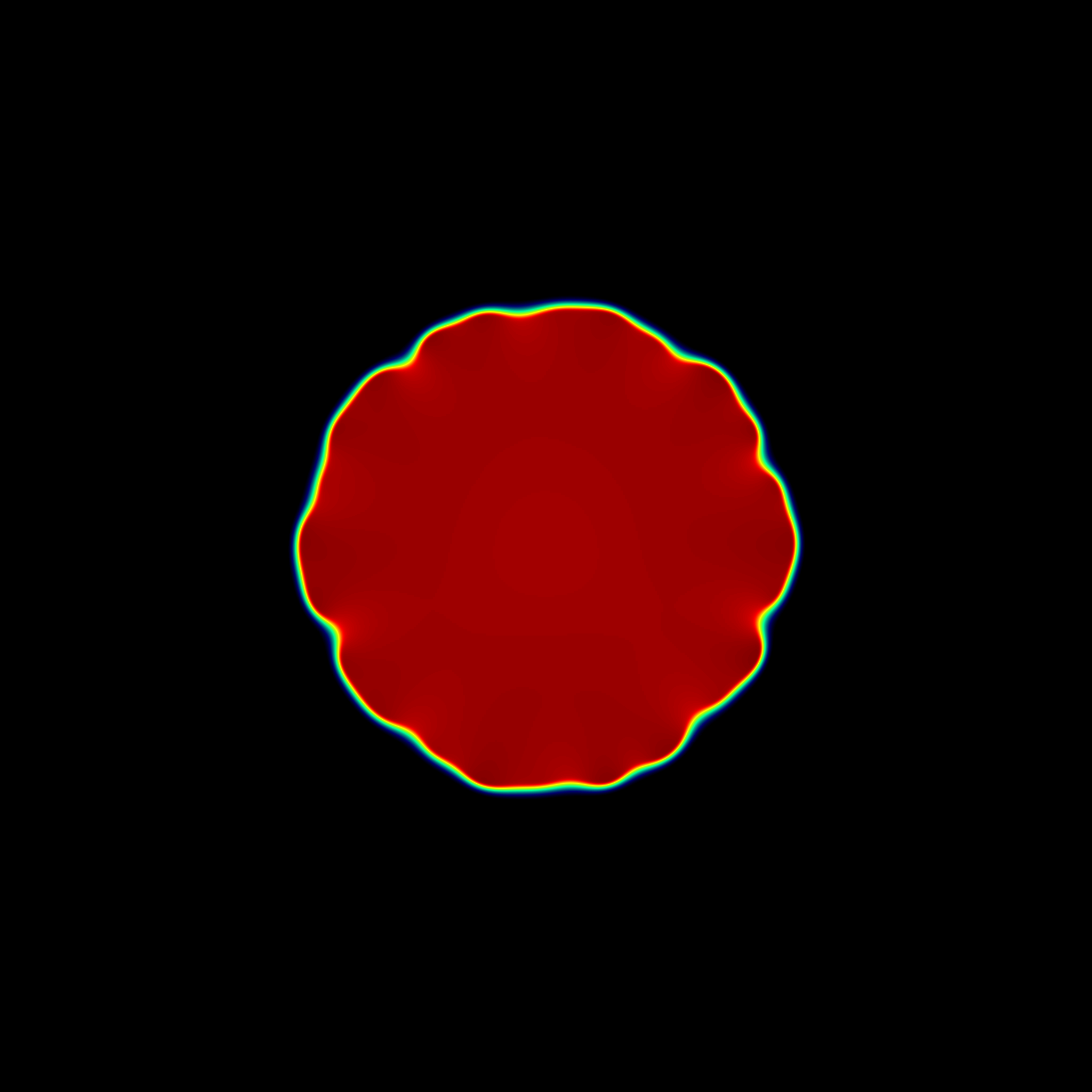}
  \caption{Samples of premixed flame kernel at final time.}
  \label{fig:kernels}
\end{figure}

\begin{figure}
  \centering
  \includegraphics[width=.5\textwidth]{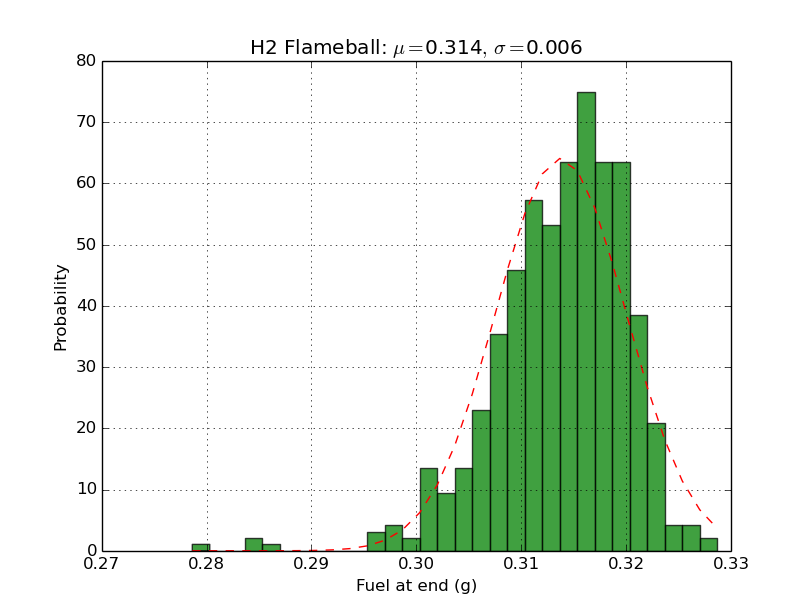}
  \caption{Histogram of products at final time for premixed flame kernel}
  \label{fig:premixed_hist}
\end{figure}

Here we select 576 samples from the posterior, at uniform intervals from the second half of the
Markov chain -- i.e., the values of the 64 walkers at step numbers $7500 + 7500 \cdot n$, for $n\in (0:8)$.
We then considered two additional computational experiments: one designed to be similar to the
experiments used to calibrate the mechanism and the second designed to explore a different
combustion regime.
For the first case, we consider the propagation of a perturbed premixed hydrogen flame kernel at lean conditions.
Specifically, we consider an open domain filled with premixed air and H$_2$ at 20~atm and
stoichiometry $\phi = 0.4$. A perturbed flame kernel is initialized at the
center of the domain and allowed to propagate outward for 130~ms.
Because this flame is thermodiffusively unstable, the detailed morphology of any finite-time
realization of the flame is expected to be strongly dependent on the initial data.  However,
the mean propagation of the integrated flame surface is less sensitive.
In Figure~\ref{fig:kernels} we show the final configuration for a set of three such flames,
indicating the typical range of variability observed over the 576 samples of the posterior.
Figure~\ref{fig:premixed_hist} shows a histogram of the integrated mass of H$_2$O at the final time
(this is a surrogate measure of the mean flame propagation throughout the domain).  The results 
indicate a fairly narrow overall distribution and, given the inherent instability of this
configuration, suggests that the mechanism has been reasonably well characterized for 
problems in this regime.

\begin{figure}
  \centering
  \includegraphics[width=.5\textwidth]{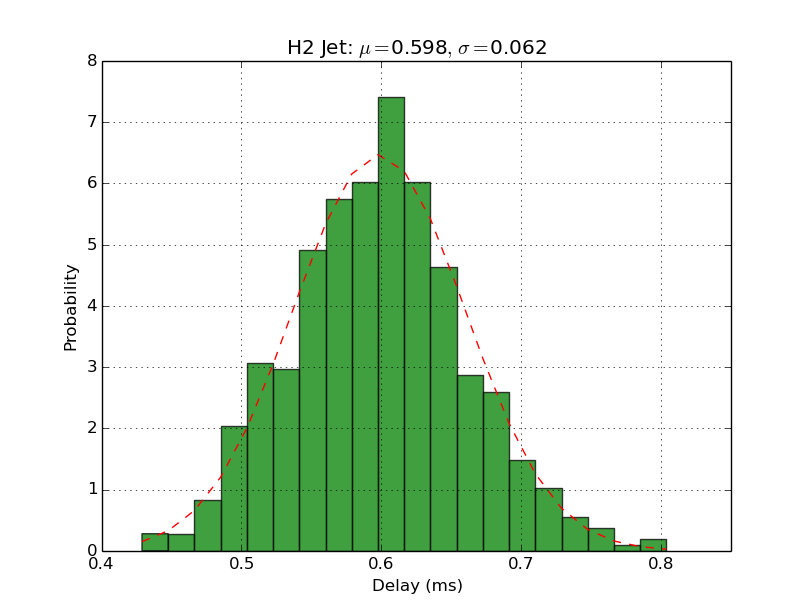}
  \caption{History of ignition delay for hydrogen jet injected into heated chamber.}
  \label{fig:jet_ign}
\end{figure}

\begin{figure}
  \centering
  \includegraphics[width=.16\textwidth]{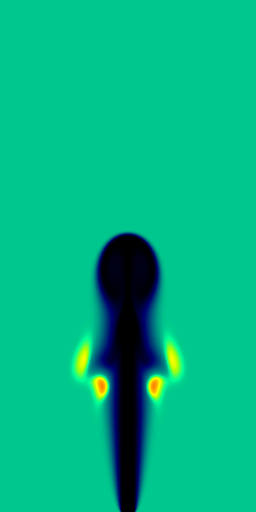}
  \includegraphics[width=.16\textwidth]{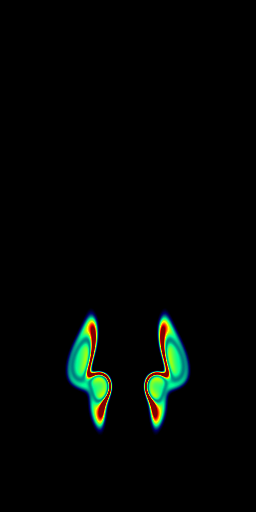}
  \includegraphics[width=.16\textwidth]{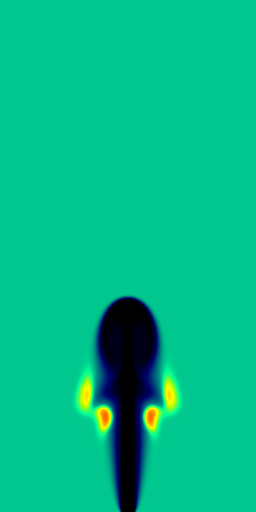}
  \includegraphics[width=.16\textwidth]{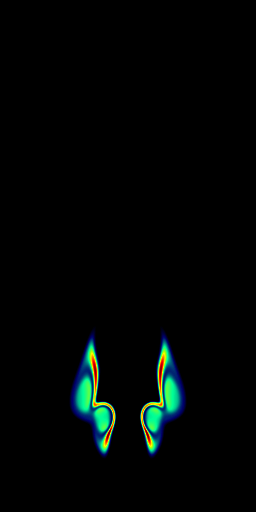}
  \includegraphics[width=.16\textwidth]{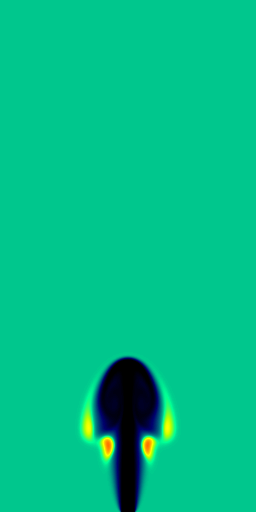}
  \includegraphics[width=.16\textwidth]{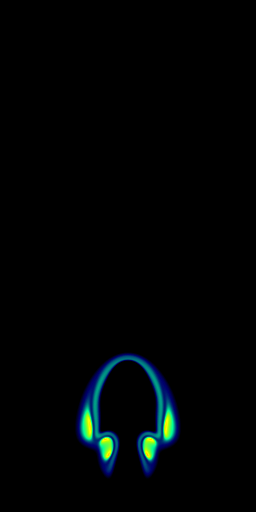}
  \caption{Samples of H$_2$ jet at ignition.  These correspond temperature and H$_2$O$_2$ mole fraction
for representative cases with long (78.8~msec), medium (60.1~msec), and short (42.9~msec) ignition times,
relative to the variation observed over the sample set.}
  \label{fig:jet_images_ign}
\end{figure}

\begin{figure}
  \centering
  \includegraphics[width=.16\textwidth]{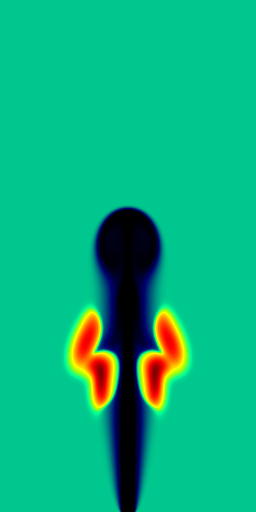}
  \includegraphics[width=.16\textwidth]{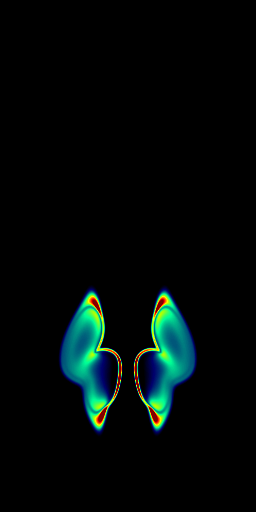}
  \includegraphics[width=.16\textwidth]{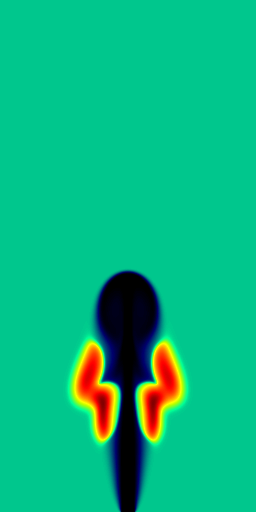}
  \includegraphics[width=.16\textwidth]{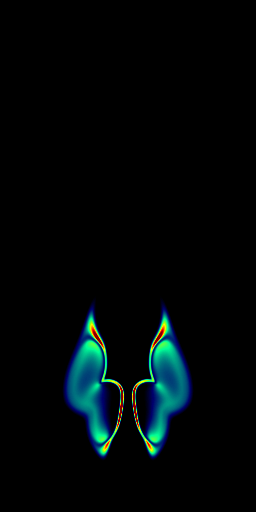}
  \includegraphics[width=.16\textwidth]{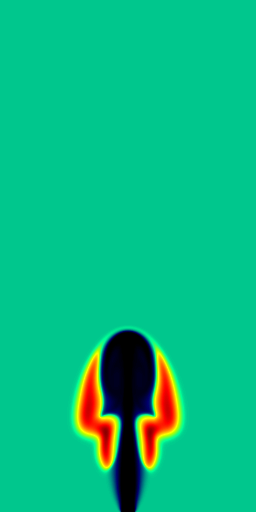}
  \includegraphics[width=.16\textwidth]{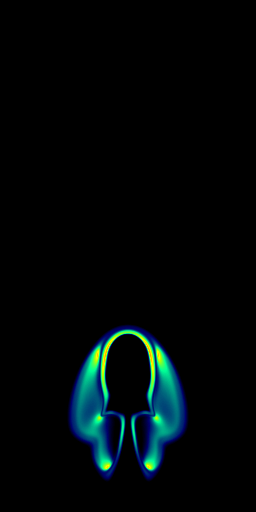}
  \caption{Temperature and H$_2$O$_2$ mole fraction of the H$_2$ jets shown in
    Figure~\ref{fig:jet_images_ign}, shown here 4~msec after ignition.}
  \label{fig:jet_images_later}
\end{figure}

For the second example, we inject a cold H$_2$ jet ($T=400$K) into a 1~atm heated chamber filled with air at
$T=1000$K.
For these simulations we are interested in the time it takes for the jet to ignite, which we define to
be the time at which the maximum temperature in the domain first exceeds $T=1850$K.
Figure \ref{fig:jet_ign} presents a histogram of the ignition time for the 576 posterior samples
discussed above.
In this case, we see substantial variability with 
a standard deviation of more than 10\% of the mean.
This variability is reflected in the details of the ignition process.
In Figure \ref{fig:jet_images_ign} we show $T$ and the mole fraction of H$_2$O$_2$
for a slow, an intermediate and and a fast ignition case from the simulated ensemble.
Figures \ref{fig:jet_images_later} shows the same quantities
four milliseconds later as
the flames continue to ignite.
In this case the images show considerable variability.  This is not surprising 
given that the bulk of the experimental data used to calibrate the mechanism
represents a different combustion regime.
More experimental data, specifically representative of this combustion regime, would likely 
reduce the variability in the predictions for this experiment.

\section{Conclusion}

We have performed a first principles Bayesian analysis of a chemical kinetic mechanism 
for hydrogen.
We have calibrated simultaneously 31 parameters effecting third body interactions and 
the $\mbox{H}_2\mbox{O}_2$
an $\mbox{H}\mbox{O}_2$
pathways that are thought to be important at high pressure.
We used data from 91 relevant experiments, 77 of which involve steady premixed flames.
For this study we assumed broad priors to isolate impact of the data on the parameter distribution.
The posterior distribution of the 31 parameters was constructed.
We generated samples from the posterior using an advanced Markov chain Monte Carlo (MCMC) algorithm.
Each evaluation of the posterior in the MCMC involved accurate simulation of all
91 experiments. 

Posterior sampling is significantly harder when using broad non-informative priors.
We felt it was important to use broad priors in order for the results to not be determined
by details of the prior, such as prior bounds on specific parameters.
Sampling is also harder in higher dimensions.
It is much harder to sample the present 31 dimensional posterior distribution than a posterior distribution involving 
significanly fewer unknown fitting parameters.
Nevertheless, we felt it was important to vary as many of the parameters
in the H$_2$O$_2$ and HO$_2$ submechanism as was feasible.
Letting parameters vary simultaneously gives a more accurate view of posterior parameter 
uncertainties.

The posterior samples allowed us to do first principles Bayesian uncertainty propagation
to two distinct complex unsteady flames: (1) an expanding lean premixed flame at high pressure, and
(2) ignition of a cold jet of fuel entering a bath of oxidizer.
For the first case, the predictions of the posterior samples were similar, as expected given that
premixed flames were well-represented amongst the calibration experiments.
In the second case, the results show considerable variability.
Distinct parameter sets consistent with all available calibration data lead to dramatically different
predictions of ignition time.  
This type of uncertainty propagation gives
more useful information about the uncertainty in parameter estimation than individual
error bars, confidence intervals, or even a posterior covariance matrix.
The posterior distribution is not well approximated by a multivariate normal.
At present we are unaware of shortcuts (response surface modeling, polynomial chaos expansions,
etc.) that would allow us to learn posterior uncertainties with this level of fidelity.
That said, it is possible that reduced order modeling of the experiments or the log likelihood 
function, together with multi-fidelity Monte Carlo techniques, could speed the sampling process.

An important area for further research is understanding why simultaneous fitting of many 
parameters is so difficult.
The MCMC sampler we are using has proven able to sample distributions in much higher dimensions.
But the present mechanism calibration problem seems significantly harder.
Flame physics and chemistry are highly non-linear and the log-likelihood surface is non-convex.
We were unable to identify simple structures such as multiple local probability peaks corresponding
to locally optimal fits.
It seems possible to us that the log-likelihood surface has a complex landscape.
We believe that a better understanding of complex log-likelihood surfaces would allow faster
samplers.
Slow sampling was the primary bottleneck in the present work.

Another useful extension of the present work would be in first principles experimental design,
also based on Bayesian uncertainty propagation.
One may be able to determine what new experiments would most reduce the uncertainty in a specific
target application or a target parameter.

It may be inappropriate to model experimental measurement errors as Gaussian.
The Gaussian distribution makes large errors less likely than they may be in practical experiments.
This may distort the results by making the calibrations too sensitive to outlier experiments.
It is not clear whether a fatter tailed distribution than the Gaussian in (\ref{eq:post})
would lead to significantly different results.

\section*{Acknowledgements}

The work at LBNL was supported by the U.S. Department of Energy, Office of Science,
Office of Advanced Scientific Computing Research, Applied Mathematics program under contract number DE-AC02005CH11231. Part of the simulations were performed using resources of the National Energy Research Scientific Computing Center (NERSC), a DOE Office of Science User Facility supported by the Office of Science of the U.S. Department of Energy under Contract No. DE-AC02-05CH11231.

\section*{References}

\bibliographystyle{unsrt}    
\bibliography{hydrogen_oxidation,chemistry_uq}

\begin{thebibliography}{10}

\bibitem{grimech}
Gregory~P. Smith, David~M. Golden, Michael Frenklach, Nigel~W. Moriarty, Boris
  Eiteneer, Mikhail Goldenberg, C.~Thomas Bowman, Ronald~K. Hanson, Soonho
  Song, Jr. William C.~Gardiner, Vitali~V. Lissianski, and Zhiwei Qin.
\newblock {GRI-Mech} 3.0.
\newblock {http://www.me.berkeley.edu/gri\_mech/}.

\bibitem{FrenklachWR92}
Michael Frenklach, Hai Wang, and Martin~J. Rabinowitz.
\newblock Optimization and analysis of large chemical kinetic mechanisms using
  the solution mapping method---combustion of methane.
\newblock {\em Progress in Energy and Combustion Science}, 18(1):47--73, 1992.

\bibitem{DavisJWE05}
Scott~G Davis, Ameya~V Joshi, Hai Wang, and Fokion Egolfopoulos.
\newblock An optimized kinetic model of {H2/CO} combustion.
\newblock {\em Proceedings of the Combustion Institute}, 30(1):1283--1292,
  2005.

\bibitem{LiYWL15}
Xiaoyu Li, Xiaoqing You, Fujia Wu, and Chung~K Law.
\newblock Uncertainty analysis of the kinetic model prediction for
  high-pressure {H2/CO} combustion.
\newblock {\em Proceedings of the Combustion Institute}, 35(1):617--624, 2015.

\bibitem{Stuart:2010}
A.~M. Stuart.
\newblock Inverse problems: A bayesian perspective.
\newblock {\em Acta Numerica}, 19:451--559, 2010.

\bibitem{BramanOR13}
Kalen Braman, Todd~A Oliver, and Venkat Raman.
\newblock Bayesian analysis of syngas chemistry models.
\newblock {\em Combustion Theory and Modelling}, 17(5):858--887, 2013.

\bibitem{emcee_hammer}
D.~{Foreman-Mackey}, D.~W. {Hogg}, D.~{Lang}, and J.~{Goodman}.
\newblock {emcee: {T}he {MCMC} {H}ammer}.
\newblock {\em Publications of the Astronomical Society of the Pacific},
  125:306--312, March 2013.

\bibitem{PREMIX98}
R.~J. Kee, J.~F. Grcar, M.~D. Smooke, and J.~A. Miller.
\newblock {PREMIX}: A {FORTRAN} program for modeling steady, laminar,
  one-dimensional premixed flames.
\newblock Technical Report SAND85-8240, Sandia National Laboratories,
  Livermore, 1983.

\bibitem{MuellerYD98}
Mark~A Mueller, Richard~A Yetter, and Frederick~L Dryer.
\newblock Measurement of the rate constant for {H+ O 2+ M}$\rightarrow${ HO 2+
  M (M= N 2, Ar)} using kinetic modeling of the high-pressure {H2/O2/NOx}
  reaction.
\newblock In {\em Symposium (International) on Combustion}, volume~27, pages
  177--184. Elsevier, 1998.

\bibitem{AshmanH98}
P.~J. Ashman and B.~S. Haynes.
\newblock Rate coefficient of {H+O$_2$ + M $\rightarrow$ HO$_2$ + M $\;$(M =
  H$_2$O, N$_2$, Ar, CO$_2$)}.
\newblock In {\em Twenty-Seventh Symposium (International) on Combustion}. The
  Combustion Institute, 1998.

\bibitem{YetterDR91}
R.~A. Yetter, F.L. Drywer, and H.~Rabitz.
\newblock {\em Combustion Science and Technology}, 79:97, 1991.

\bibitem{MuellerKYD99}
MA~Mueller, TJ~Kim, RA~Yetter, and FL~Dryer.
\newblock Flow reactor studies and kinetic modeling of the {H2/O2} reaction.
\newblock {\em International Journal of Chemical Kinetics}, 31(2):113--125,
  1999.

\bibitem{LiZKD04}
Juan Li, Zhenwei Zhao, Andrei Kazakov, and Frederick~L Dryer.
\newblock An updated comprehensive kinetic model of hydrogen combustion.
\newblock {\em International journal of chemical kinetics}, 36(10):566--575,
  2004.

\bibitem{ConaireCSPW04}
Marcus {\'O}~Conaire, Henry~J Curran, John~M Simmie, William~J Pitz, and
  Charles~K Westbrook.
\newblock A comprehensive modeling study of hydrogen oxidation.
\newblock {\em International Journal of Chemical Kinetics}, 36(11):603--622,
  2004.

\bibitem{BurkeCDJ10}
Michael~P Burke, Marcos Chaos, Frederick~L Dryer, and Yiguang Ju.
\newblock Negative pressure dependence of mass burning rates of
  {H2/CO/O2/diluent} flames at low flame temperatures.
\newblock {\em Combustion and Flame}, 157(4):618--631, 2010.

\bibitem{BurkeCJDK12}
Michael~P Burke, Marcos Chaos, Yiguang Ju, Frederick~L Dryer, and Stephen~J
  Klippenstein.
\newblock Comprehensive {H2/O2} kinetic model for high-pressure combustion.
\newblock {\em International Journal of Chemical Kinetics}, 44(7):444--474,
  2012.

\bibitem{Keromnes_Many_Curran13}
Alan K{\'e}romn{\`e}s, Wayne~K Metcalfe, Karl~A Heufer, Nicola Donohoe,
  Apurba~K Das, Chih-Jen Sung, J{\"u}rgen Herzler, Clemens Naumann, Peter
  Griebel, Olivier Mathieu, et~al.
\newblock An experimental and detailed chemical kinetic modeling study of
  hydrogen and syngas mixture oxidation at elevated pressures.
\newblock {\em Combustion and Flame}, 160(6):995--1011, 2013.

\bibitem{Konnov08}
Alexander~A Konnov.
\newblock Remaining uncertainties in the kinetic mechanism of hydrogen
  combustion.
\newblock {\em Combustion and flame}, 152(4):507--528, 2008.

\bibitem{LiZKCDS07}
Juan Li, Zhenwei Zhao, Andrei Kazakov, Marcos Chaos, Frederick~L Dryer, and
  James~J Scire.
\newblock A comprehensive kinetic mechanism for {CO, CH2O, and CH3OH}
  combustion.
\newblock {\em International Journal of Chemical Kinetics}, 39(3):109--136,
  2007.

\bibitem{YouPF11}
Xiaoqing You, Andrew Packard, and Michael Frenklach.
\newblock Process informatics tools for predictive modeling: Hydrogen
  combustion.
\newblock {\em International Journal of Chemical Kinetics}, 44(2):101--116,
  2012.

\bibitem{CHEMKINIII96}
R.~J. Kee, F.~M. Ruply, E.~Meeks, and J.~A. Miller.
\newblock {CHEMKIN-III: A FORTRAN} chemical kinetics package for the analysis
  of gas-phase chemical and plasma kinetics.
\newblock Technical Report SAND96-8216, 1996.

\bibitem{vode}
P.~N. Brown, G.~D. Byrne, and A.~C. Hindmarsh.
\newblock {VODE}: A variable coefficient ode solver.
\newblock {\em SIAM J. Sci. Stat. Comput.}, 10:1038--1051, 1989.

\bibitem{GoodmanWeare:2010}
J.~Goodman and J.~Weare.
\newblock Ensemble samplers with affine invariance.
\newblock {\em Comm. App. Math. and Comp. Sci.}, pages 65--80, 2010.

\bibitem{Sokal:1989}
Jonathan Goodman and Alan~D. Sokal.
\newblock Multigrid {Monte Carlo} method. {C}onceptual foundations.
\newblock {\em Phys. Rev. D}, 40:2035--2071, Sep 1989.

\bibitem{BurkeDJ11}
Michael~P Burke, Frederick~L Dryer, and Yiguang Ju.
\newblock Assessment of kinetic modeling for lean {H2/CH4/O2/diluent} flames at
  high pressures.
\newblock {\em Proceedings of the Combustion Institute}, 33(1):905--912, 2011.

\end{thebibliography}

\end{document}